\documentclass[12pt]{amsart}

\usepackage{amssymb}
\usepackage{ifthen}
\usepackage{mathrsfs}
\usepackage{pict2e}
\usepackage{xargs}
\usepackage{xspace}

\newcommand{\definedterm}[1]{\emph{#1}}

\newcommand{\branches}[1]{[#1]}
\newcommand{\calF}{\mathscr{F}}
\newcommand{\calI}{\mathcal{I}}
\newcommand{\calK}{\mathcal{K}}
\newcommand{\calN}{\mathcal{N}}
\newcommand{\calP}{\mathcal{P}}
\newcommand{\Cantorspace}[1][]{\ifthenelse{\equal{#1}{}}{\functions{\N}{2}}{\functions{#1}{2}}}
\newcommand{\CantorCantorspace}{\functions{\N}{(\Cantorspace)}}
\newcommand{\Cantortree}{\functions{<\N}{2}}
\newcommandx{\cardinality}[2][2 =]{\ifthenelse{\equal{#2}{}}{|#1|}{|#1|_{#2}}}
\newcommand{\cardinalityideal}[1]{\calK_{#1}}
\newcommand{\characteristicfunction}[1]{\chi_{#1}}
\newcommand{\closure}[1]{\mathrm{cl}(#1)}
\newcommand{\composition}{\circ}
\newcommandx{\concatenation}[2][1 = undefined, 2 = undefined]{
  \ifthenelse{\equal{#1}{undefined}}{{}\smallfrown}{
    \ifthenelse{\equal{#2}{undefined}}{\bigoplus #1}{\bigoplus_{#1} #2}
  } 
}

\newcommand{\emptystring}{\emptyset}
\newcommand{\Eone}{\mathbb{E}_1}
\newcommand{\equivalenceclass}[2]{[#1]_{#2}}
\newcommand{\extendedby}{\sqsubseteq}
\newcommand{\extensions}[1]{\calN_{#1}}
\newcommand{\Ezero}{\mathbb{E}_0}
\newcommand{\from}{\colon}
\newcommand{\functions}[2]{#2^{#1}}
\newcommand{\image}[2]{#1(#2)}
\newcommand{\inducedequivalencerelation}[1]{E_{#1}}
\newcommand{\inducedgraph}[2]{#1_{#2}}
\newcommandx{\intersection}[2][1 = undefined, 2 = undefined]{
  \ifthenelse{\equal{#1}{undefined}}{\cap}{
    \ifthenelse{\equal{#2}{undefined}}{\bigcap #1}{\bigcap_{#1} #2}
  }
}
\newcommand{\mathand}{\text{ and }}
\newcommand{\N}{\mathbb{N}}
\newcommand{\pair}[2]{(#1, #2)}
\newcommandx{\Piclass}[2][1=, 2=]{\mathbf{\Pi}^{\pmb{#1}}_{\pmb{#2}}}
\newcommand{\powerset}[1]{\calP(#1)}
\newcommand{\preimage}[2]{#1^{-1}(#2)}
\newcommand{\projection}[1]{\mathrm{proj}_{#1}}
\newcommand{\pushforward}[2]{#1_* #2}
\newcommand{\saturation}[2]{[#1]_{#2}}
\newcommand{\scrM}{\mathscr{M}}
\newcommand{\scrN}{\mathscr{N}}
\newcommandx{\sequence}[2][2 = undefined]{
  \ifthenelse{\equal{#2}{undefined}}{(#1)}{
    (#1)_{#2}
  }
}
\newcommandx{\set}[2][2 = undefined]{
  \ifthenelse{\equal{#2}{undefined}}{\{ #1 \}}{
    \{ #1 \suchthat #2 \}
  }
}
\newcommandx{\sets}[4][3 = undefined, 4 = undefined]{
  \ifthenelse{\equal{#3}{undefined}}{[#2]^{#1}}{
    \ifthenelse{\equal{#4}{undefined}}{[#2]^{#1}_{#3}}{[#2]^{#1}_{#3, #4}}
  }
}
\newcommandx{\Sigmaclass}[2][1=, 2=]{\mathbf{\Sigma}^{\pmb{#1}}_{\pmb{#2}}}
\newcommand{\strictlyextendedby}{\sqsubset}
\newcommand{\suchthat}{\mid}
\renewcommand{\restriction}[2]{#1 \upharpoonright #2}
\newcommand{\singletonsequence}[1]{\sequence{#1}}
\newcommand{\support}[1]{\mathrm{supp}(#1)}
\newcommand{\symmetricdifference}{\mathbin{\triangle}}
\newcommand{\textexponent}[2]{$#1^{\mathrm{#2}}$}
\newcommandx{\union}[2][1 = undefined, 2 = undefined]{
  \ifthenelse{\equal{#1}{undefined}}{\cup}{
    \ifthenelse{\equal{#2}{undefined}}{\bigcup #1}{\bigcup_{#1} #2}
  }
}
\newcommand{\verticalsection}[2]{#1_{#2}}

\newcommand{\Borel}{Bor\-el\xspace}
\newcommand{\Dougherty}{Dough\-er\-ty\xspace}
\newcommand{\Feldman}{Feld\-man\xspace}
\newcommand{\Harrington}{Har\-ring\-ton\xspace}
\newcommand{\Hjorth}{Hjorth\xspace}
\newcommand{\Jackson}{Jack\-son\xspace}
\newcommand{\Kechris}{Kech\-ris\xspace}
\newcommand{\Louveau}{Lou\-veau\xspace}
\newcommand{\Lusin}{Lu\-sin\xspace}
\newcommand{\Moore}{Moore\xspace}
\newcommand{\Novikov}{No\-vik\-ov\xspace}
\newcommand{\Polish}{Po\-lish\xspace}
\newcommand{\Souslin}{Sous\-lin\xspace}
\newcommand{\Wadge}{Wadge\xspace}

\newenvironment{lemmaproof}{
   
  \begin{proof}
}{\end{proof}}

\newenvironment{propositionproof}{
   
  \begin{proof}
}{\end{proof}}

\newenvironment{theoremproof}{
   
  \begin{proof}
}{\end{proof}}

\newtheorem{introtheorem}{Theorem}
\newtheorem{lemma}{Lemma}[section]
\newtheorem{proposition}[lemma]{Proposition}
\newtheorem{theorem}[lemma]{Theorem}

\theoremstyle{definition}
\newtheorem*{acknowledgments}{Acknowledgments}

\begin{document}

\baselineskip=13.97pt

\author[J.D. Clemens]{John D. Clemens}

\address{
  John D. Clemens \\
  Southern Illinois University \\
  Mathematics Department \\
  Neckers \\
  1245 Lincoln Drive \\
  Mailstop 4408 \\
  Carbondale, IL 62901 \\
  USA
 }

\email{clemens@siu.edu}

\urladdr{
  http://www.math.siu.edu/faculty-staff/faculty/clemens.php
}

\author[D. Lecomte]{Dominique Lecomte}

\address{
  Dominique Lecomte \\
  Universit\'{e} Paris 6, Institut de Math\'{e}matiques de Jussieu, Projet Analyse Fonctionnelle \\
  Couloir 16 - 26, 4\`{e}me \'{e}tage, Case 247, 4, place Jussieu, 75 252 Paris Cedex 05, France \\
  and Universit\'{e} de Picardie, I.U.T de l'Oise, site de Creil \\
  13, all\'{e}e de la fa\"iencerie, 60 107 Creil, France
}

\email{dominique.lecomte@upmc.fr}

\urladdr{
  https://www.imj-prg.fr/~dominique.lecomte/
}

\author[B.D. Miller]{Benjamin D. Miller}

\address{
  Benjamin D. Miller \\
  Kurt G\"{o}del Research Center for Mathematical Logic \\
  Universit\"{a}t Wien \\
  W\"{a}hringer Stra{\ss}e 25 \\
  1090 Wien \\
  Austria \\ 
  and Institut f\"{u}r mathematische Logik und Grundlagenforschung \\
  Fachbereich Mathematik und Informatik \\
  Universit\"{a}t M\"{u}nster \\
  Einsteinstra{\ss}e 62 \\
  48149 M\"{u}nster \\
  Germany
 }

\email{glimmeffros@gmail.com}

\urladdr{
  http://wwwmath.uni-muenster.de/u/ben.miller
}

\thanks{The authors were supported in part by SFB Grant 878.}

  

\title
  [Dichotomies and essential complexity]
  {Dichotomy theorems for families of non-cofinal essential complexity}

\begin{abstract}
  We prove that for every \Borel equivalence relation $E$, either $E$ is \Borel reducible to 
  $\Ezero$, or the family of \Borel equivalence relations incompatible with $E$ has cofinal
  essential complexity. It follows that if $F$ is a \Borel equivalence relation and $\calF$ is a 
  family of \Borel equivalence relations of non-cofinal essential complexity which together
  satisfy the dichotomy that for every \Borel equivalence relation $E$, either $E \in \calF$ or
  $F$ is \Borel reducible to $E$, then $\calF$ consists solely of smooth equivalence relations,
  thus the dichotomy is equivalent to a known theorem.
\end{abstract}

\maketitle

\section*{Introduction}

A \definedterm{reduction} of an equivalence relation $E$ on a set $X$ to an equivalence relation
$F$ on a set $Y$ is a function $\pi \from X \to Y$ with the property that $\forall x_1, x_2 \in X
\ ( x_1 \mathrel{E} x_2 \iff \pi(x_1) \mathrel{F} \pi(x_2) )$. A topological space is \definedterm{\Polish}
if it is second countable and completely metrizable, a subset of such a space is \definedterm
{\Borel} if it is in the $\sigma$-algebra generated by the underlying topology, and a function between
such spaces is \definedterm{\Borel} if pre-images of open sets are \Borel. Over the last few 
decades, the study of \Borel reducibility of \Borel equivalence relations on \Polish spaces has 
emerged as a central theme in descriptive set theory.

The early development of this area was dominated by dichotomy theorems. There are several
trivial ones, such as the fact that if $n$ is a natural number, then for every \Borel equivalence 
relation $E$ on a \Polish space, either $E$ is \Borel reducible to equality on $n$, or equality on 
$n + 1$ is \Borel reducible to $E$. Similarly, either there is a natural number $n$ for which $E$ is 
\Borel reducible to equality on $n$, or equality on $\N$ is \Borel reducible to $E$.

There are also non-trivial results of this form. By \cite{Silver}, either $E$ is \Borel reducible to
equality on $\N$, or equality on $\Cantorspace$ is \Borel reducible to $E$. And by \cite
[Theorem 1.1]{HarringtonKechrisLouveau}, either $E$ is \Borel reducible to equality on 
$\Cantorspace$, or $\Ezero$ is \Borel reducible to $E$, where $\Ezero$ is the relation on
$\Cantorspace$ given by $x \mathrel{\Ezero} y \iff \exists n \in \N \forall m \ge n \ x(m) = y(m)$.

Whereas the results we have mentioned thus far concern the global structure of the \Borel
reducibility hierarchy, \cite[Theorem 1]{KechrisLouveau} yields a local dichotomy of this form.
Namely, that for every \Borel equivalence relation $E$ on a \Polish space which is \Borel reducible 
to $\Eone$, either $E$ is \Borel reducible to $\Ezero$, or $\Eone$ is \Borel reducible to $E$, where 
$\Eone$ is the relation on $\CantorCantorspace$ given by $x \mathrel{\Eone} y \iff \exists n \in \N 
\forall m \ge n \ x(m) = y(m)$.

At first glance, one might hope the assumption that $E$ is \Borel reducible to $\Eone$ could be
eliminated, thereby yielding a new global dichotomy theorem. Unfortunately, \cite[Theorem 2]
{KechrisLouveau} ensures that if $E$ is not \Borel reducible to $\Ezero$, then there is a \Borel 
equivalence relation with which it is incomparable under \Borel reducibility. It follows that only the 
pairs $\pair{F}{F'}$ discussed thus far (up to \Borel bi-reducibility) satisfy both (1) there is a \Borel
reduction of $F$ to $F'$ but not vice versa, and (2) for every \Borel equivalence relation $E$ on a
\Polish space, either $E$ is \Borel reducible to $F$, or $F'$ is \Borel reducible to $E$.

As the latter result rules out further global dichotomies of the sort discussed thus far, it is
interesting to note that its proof hinges upon the previously mentioned local dichotomy, as well as
\Harrington's unpublished theorem that the family of orbit equivalence relations induced by
\Borel actions of \Polish groups on \Polish spaces is unbounded in the \Borel reducibility
hierarchy. Here we utilize strengthenings of these results to provide a substantially stronger 
anti-dichotomy theorem.

Given a property $P$ of \Borel equivalence relations, we say that a \Borel equivalence relation
is \definedterm{essentially $P$} if it is \Borel reducible to a \Borel equivalence relation on a \Polish
space with the given property. A \definedterm{\Wadge reduction} of a set $A \subseteq X$ to a set
$B \subseteq Y$ is a continuous function $\pi \from X \to Y$ such that $\forall x \in X \ ( x \in A \iff 
\pi(x) \in B )$. We say that a \Borel equivalence relation $E$ has \definedterm{essential complexity 
at least the complexity} of a set $B \subseteq \Cantorspace$ if $B$ is \Wadge reducible to every
\Borel equivalence relation to which $E$ is \Borel reducible. We say that a family $\calF$ of
\Borel equivalence relations has \definedterm{cofinal essential complexity} if for every 
\Borel set $B \subseteq \Cantorspace$, there is a \Borel equivalence relation $E \in \calF$ with
essential complexity at least the complexity of $B$.

Much as in \cite{KechrisLouveau}, we obtain our anti-dichotomy theorem as a consequence of a
result yielding the existence of incomparable \Borel equivalence relations, albeit one considerably 
stronger than that given there.

\begin{introtheorem} \label{introduction:main}
  Suppose that $X$ is a \Polish space and $E$ is a \Borel equivalence relation on $X$. Then exactly
  one of the following holds:
  \begin{enumerate}
    \item The relation $E$ is \Borel reducible to $\Ezero$.
    \item The family of \Borel equivalence relations on \Polish spaces which are incomparable with
      $E$ under \Borel reducibility has cofinal essential complexity.
  \end{enumerate}
\end{introtheorem}

We say that a family $\calF$ of \Borel equivalence relations on \Polish spaces is \definedterm
{dichotomical} if there is a \definedterm{minimum} \Borel equivalence relation $F$ on a \Polish 
space which is not in $\calF$, meaning that whenever $E$ is a \Borel equivalence relation on a
\Polish space, either $E \in \calF$ or there is a \Borel reduction of $F$ to $E$. The following 
consequence of Theorem \ref{introduction:main} implies that the only such families are those 
associated with the dichotomies mentioned earlier.

\begin{introtheorem}
  Suppose that $\calF$ is a dichotomical class of \Borel equivalence relations on \Polish spaces of 
  non-cofinal essential complexity. Then every equivalence relation in $\calF$ is smooth.
\end{introtheorem}

In \S\ref{preliminaries}, we briefly review the preliminaries needed throughout the paper. In
\S\ref{hyperfiniteoncountable}, we introduce a property of graphs $G$ ensuring that if a \Borel 
equivalence relation $E$ on a \Polish space, whose classes are all countable, is \Borel reducible
to the equivalence relation generated by $G$, then it is \Borel reducible to $\Ezero$. In \S\ref{ideals}, 
we introduce a family of ideals on $\N \times \N$, and show that they are cofinal under \Wadge 
reducibility. In \S\ref{trees}, we introduce a family of trees on $\N \times \N$, and show that the 
graphs determined by their branches interact nicely with equivalence relations induced by ideals.
In \S\ref{density}, we consider a subfamily of these trees satisfying an appropriate density
condition, and show that the graphs determined by their branches interact particularly nicely with 
equivalence relations induced by the ideals introduced earlier. And in \S\ref{antibasis}, we establish
our primary results.

\section{Preliminaries} \label{preliminaries}

Two sets $M, N \subseteq \N$ are \definedterm{almost disjoint} if $\cardinality{M \intersection N} < 
\aleph_0$.

\begin{proposition} \label{preliminaries:almostdisjoint}
  There is a continuous injection $\pi \from \Cantorspace \to \powerset{\N}$ into a family of pairwise 
  almost disjoint infinite sets.
\end{proposition}

\begin{propositionproof}
  It is sufficient to observe that the function $\pi \from \Cantorspace \to \powerset{\Cantortree}$, given
  by $\pi(c) = \set{\restriction{c}{n}}[n \in \N]$, is a continuous injection into a family of pairwise almost
  disjoint infinite sets.
\end{propositionproof}

A set is \definedterm{comeager} if it contains an intersection of countably many dense open sets.

\begin{proposition} \label{preliminaries:continuous}
  Suppose that $X$ and $Y$ are \Polish spaces and the map $\pi \from X \to Y$ is \Borel . Then
  there is a comeager set $C \subseteq X$ on which $\pi$ is continuous.
\end{proposition}

\begin{propositionproof}
  See, for example, \cite[Theorem 8.38]{Kechris}.
\end{propositionproof}

A subset of a \Polish space is \definedterm{analytic} if it is the continuous image of a \Borel subset of
a \Polish space. A subset of a \Polish space is \definedterm{co-analytic} if its complement is analytic.

\begin{theorem}[\Souslin] \label{preliminaries:bianalytic}
  Suppose that $X$ is a \Polish space and $B \subseteq X$. Then $B$ is \Borel if and only if $B$ is
  both analytic and co-analytic.
\end{theorem}

\begin{theoremproof}
  See, for example, \cite[Theorem 14.11]{Kechris}.
\end{theoremproof}

The \definedterm{projection} from $X \times Y$ to $X$ is given by $\projection{X}(x, y) = x$. A 
\definedterm{partial uniformization} of a set $R \subseteq X \times Y$ is a function whose graph 
is contained in $R$. A \definedterm{uniformization} of a set $R \subseteq X \times Y$ is a partial 
uniformization of $R$ whose domain is $\image{\projection{X}}{R}$.

\begin{theorem}[\Lusin-\Novikov] \label{preliminaries:LusinNovikov}
  Suppose that $X$ and $Y$ are \Polish spaces and $R \subseteq X \times Y$ is a \Borel set 
  whose vertical sections are all countable. Then $\image{\projection{X}}{R}$ is \Borel, and 
  $R$ is a countable union of \Borel uniformizations.
\end{theorem}

\begin{theoremproof}
  See, for example, \cite[Theorem 18.10]{Kechris}.
\end{theoremproof}

For each $x \in X$, the \definedterm{\textexponent{x}{th} vertical section} of a set $R \subseteq
X \times Y$ is given by $\verticalsection{R}{x} = \set{y \in Y}[x \mathrel{R} y]$. The \definedterm
{set of unicity} of $R$ is $\set{x \in X}[\cardinality{\verticalsection{R}
{x}} = 1]$.

\begin{theorem}[\Lusin] \label{preliminaries:unicity}
  Suppose that $X$ and $Y$ are \Polish spaces and $R \subseteq X \times Y$ is \Borel. Then
  the set of unicity of $R$ is co-analytic.
\end{theorem}

\begin{theoremproof}
  See, for example, \cite[Theorem 18.11]{Kechris}.
\end{theoremproof}

A \definedterm{graph} on a set $X$ is an irreflexive, symmetric set $G \subseteq X \times X$.
Such a graph is \definedterm{locally countable} if its vertical sections are countable.
An \definedterm{edge $N$-coloring} of $G$ is a map $c \from G \to N$ with $\forall
\pair{x}{y} \in G \ c(x, y) = c(y, x)$ and $\forall \pair{x}{y}, \pair{x}{z}  \in G \ (y \neq z \implies
c(x, y) \neq c(x, z))$.

\begin{theorem}[\Feldman-\Moore] \label{preliminaries:coloring}
  Suppose that $X$ is a \Polish space and $G$ is a locally countable \Borel graph on $X$. Then there
  is a \Borel edge $\N$-coloring of $G$.
\end{theorem}

\begin{theoremproof}
  This follows from the proof of \cite[Theorem 1]{FeldmanMoore}.
\end{theoremproof}

We say that a \Borel equivalence relation is \definedterm{smooth} if it is \Borel reducible to equality
on $\Cantorspace$. An \definedterm{embedding} is an injective reduction.

\begin{theorem}[\Harrington-\Kechris-\Louveau] \label{preliminaries:Ezero}
  Suppose that $X$ is a \Polish space and $E$ is a \Borel equivalence relation on $X$. Then exactly
  one of the following holds:
  \begin{enumerate}
    \item The relation $E$ is smooth.
    \item There is a continuous embedding of $\Ezero$ into $E$.
  \end{enumerate}
\end{theorem}

\begin{theoremproof}
  See \cite[Theorem 1.1]{HarringtonKechrisLouveau}.
\end{theoremproof}

A \definedterm{partial transversal} of an equivalence relation $E$ on $X$ is a set $B \subseteq X$ 
intersecting every equivalence class of $E$ in at most one point. A \definedterm{transversal} of an 
equivalence relation $E$ on $X$ is a set $B \subseteq X$ intersecting every equivalence class of
$E$ in exactly one point.

Following the standard abuse of language, we say that an equivalence relation is
\definedterm{countable} if all of its equivalence classes are countable.

\begin{proposition} \label{preliminaries:countablesmooth}
  Suppose that $X$ is a \Polish space and $E$ is a countable \Borel equivalence relation on $X$.
  Then $E$ is smooth if and only if $X$ is the union of countably many \Borel partial transversals of 
  $E$.
\end{proposition}

\begin{propositionproof}
  This is a straightforward consequence of Theorem \ref{preliminaries:LusinNovikov}.
\end{propositionproof}

Again following the standard abuse of language, we say that an equivalence relation is
\definedterm{finite} if all of its equivalence classes are finite.

\begin{proposition} \label{preliminaries:finite}
  Suppose that $X$ is a \Polish space and $E$ is a finite \Borel equivalence relation on $X$. Then
  $E$ is smooth.
\end{proposition}

\begin{propositionproof}
  This is also a straightforward consequence of Theorem \ref{preliminaries:LusinNovikov}.
\end{propositionproof}

We say that a \Borel equivalence relation is \definedterm{hyperfinite} if it is the union of an increasing 
sequence $\sequence{F_n}[n \in \N]$ of finite \Borel subequivalence relations. By \cite[Theorem 1]
{DoughertyJacksonKechris}, a countable \Borel equivalence relation is hyperfinite if and only if it is 
\Borel reducible to $\Ezero$.

\begin{proposition}[Dougherty-Jackson-Kechris] \label{preliminaries:hyperfinite}
  Suppose that $X$ and $Y$ are \Polish spaces, $E$ and $F$ are countable \Borel equivalence
  relations on $X$ and $Y$, $E$ is \Borel reducible to $F$, and $F$ is hyperfinite. Then $E$ is
  hyperfinite.
\end{proposition}

\begin{propositionproof}
  See \cite[Proposition 5.2]{DoughertyJacksonKechris}.
\end{propositionproof}

We say that a \Borel equivalence relation is \definedterm{hypersmooth} if it is the union of an
increasing sequence $\sequence{F_n}[n \in \N]$ of smooth \Borel subequivalence relations. By \cite
[Propositions 1.1 and 1.3]{KechrisLouveau}, a \Borel equivalence relation is hypersmooth if and only
if it is \Borel reducible to $\Eone$.

\begin{theorem}[\Dougherty-\Jackson-\Kechris] \label{preliminaries:hypersmooth}
  Suppose that $X$ is a \Polish space and $E$ is a \Borel equivalence relation on $X$. Then $E$
  is hyperfinite if and only if $E$ is both countable and hypersmooth.
\end{theorem}

\begin{theoremproof}
  See, for example, \cite[Theorem 5.1]{DoughertyJacksonKechris}.
\end{theoremproof}

A property $P$ of subsets of $Y$ is \definedterm{$\Piclass[1][1]$-on-$\Sigmaclass[1][1]$} if
whenever $X$ is a \Polish space and $R \subseteq X \times Y$ is an analytic set, the
corresponding set $\set{x \in X}[\verticalsection{R}{x} \text{ satisfies } P]$ is co-analytic.

\begin{theorem} \label{preliminaries:reflection}
  Suppose that $X$ is a \Polish space, $\Phi$ is a $\Piclass[1][1]$-on-$\Sigmaclass[1][1]$ property 
  of subsets of $X$, and $A \subseteq X$ is an analytic set on which $\Phi$ holds. Then there is a 
  \Borel set $B \supseteq A$ on which $\Phi$ holds.
\end{theorem}

\begin{theoremproof}
  See, for example, \cite[Theorem 35.10]{Kechris}.
\end{theoremproof}

A \definedterm{path} through a graph $G$ is a sequence $\sequence{x_i}[i \le n]$ with the
property that $\forall i < n \ x_i \mathrel{G} x_{i+1}$. Such a path is a \definedterm{cycle} if
$n > 2$, $\sequence{x_i}[i < n]$ is injective, and $x_0 = x_n$. A graph is \definedterm{acyclic}
if it has no cycles.

The equivalence relation \definedterm{generated} by a graph $G$ on a set $X$ is the smallest
equivalence relation $\inducedequivalencerelation{G}$ on $X$ containing it. A \definedterm
{graphing} of an equivalence relation is a graph generating it. We say that a \Borel equivalence
relation is \definedterm{treeable} if it has an acyclic \Borel graphing.

\begin{theorem}[\Hjorth] \label{preliminaries:treeable}
  Suppose that $X$ is a \Polish space and $E$ is a treeable \Borel equivalence relation on $X$.
  Then the following are equivalent:
  \begin{enumerate}
    \item The relation $E$ is essentially countable.
    \item There is a \Borel set $B \subseteq X$ whose intersection with each equivalence class of
      $E$ is countable and non-empty.
  \end{enumerate}
\end{theorem}

\begin{theoremproof}
  See \cite[Theorem 6]{Hjorth}.
\end{theoremproof}

\begin{theorem} \label{preliminaries:dichotomy}
  Suppose that $X$ is a \Polish space and $E$ is a \Borel equivalence relation on $X$
  essentially generated by a \Borel subgraph of an acyclic compact graph. Then exactly one
  of the following holds:
  \begin{enumerate}
    \item The relation $E$ is essentially countable.
    \item There is a continuous embedding of $\Eone$ into $E$.
  \end{enumerate}
\end{theorem}

\begin{theoremproof}
  See \cite[Theorem 6.3 and Proposition 6.4]{ClemensLecomteMiller}.
\end{theoremproof}

A topological group is \definedterm{\Polish} if it is \Polish as a topological space.

\begin{theorem}[\Hjorth-\Kechris-\Louveau] \label{preliminaries:groupaction}
  The family of orbit equivalence relations induced by \Borel actions of \Polish groups on
  \Polish spaces has cofinal essential complexity.
\end{theorem}

\begin{theoremproof}
  See \cite[Theorem 4.1]{HjorthKechrisLouveau}.
\end{theoremproof}

\begin{theorem}[\Kechris-\Louveau] \label{preliminaries:nonreducible}
  Suppose that $X$ is a \Polish space and $E$ is the orbit equivalence relation induced by a \Borel
  action of a \Polish group on a \Polish space. Then $\Eone$ is not \Borel reducible to $E$.
\end{theorem}

\begin{theoremproof}
  See \cite[Theorem 4.2]{KechrisLouveau}.
\end{theoremproof}

\section{Partition stratifications} \label{hyperfiniteoncountable}

Suppose that $X$ is a \Polish space, $E$ is a \Borel equivalence relation on $X$, and $G$ is a
\Borel graph on $X$. We use $\inducedgraph{G}{E}$ to denote the graph on $X / E$ given by 
\begin{equation*}
  \inducedgraph{G}{E} = \set{\pair{C}{D} \in (X / E) \times (X / E)}[C \neq D \mathand 
    (C \times D) \intersection G \neq \emptyset].
\end{equation*}
We say that $G$ has \definedterm{faithful cycles} over $E$ if among all $G$-paths $\sequence{x_i}
[i \le k]$, only $G$-cycles have the property that $\sequence{\equivalenceclass{x_i}{E}}[i \le k]$ is 
a $\inducedgraph{G}{E}$-cycle.

A \definedterm{partition stratification} of $G$ is a sequence of the form $\sequence{E_n, G_n}[n \in 
\N]$, where $\sequence{E_n}[n \in \N]$ is a decreasing sequence of equivalence relations on $X$, 
each of which having only countably many classes, whose intersection is the diagonal,
$\sequence{G_n}[n \in \N]$ is an increasing sequence of \Borel graphs whose union is $G$, and 
each $G_n$ has faithful cycles on $E_n$.

\begin{proposition} \label{hyperfiniteoncountable:countablestratification}
  Suppose that $X$ is a \Polish space, $E$ is a countable \Borel equivalence relation on $X$, and
  $G$ is a \Borel graphing of $E$. Then $E$ is hyperfinite if and only if there is a partition 
  stratification of $G$.
\end{proposition}

\begin{propositionproof}
  Suppose first that $E$ is hyperfinite, and fix an increasing sequence $\sequence{F_n}[n \in \N]$ of 
  finite \Borel subequivalence relations whose union is $E$. As Propositions \ref
  {preliminaries:countablesmooth} and \ref{preliminaries:finite} ensure that spaces underlying finite 
  \Borel equivalence relations are countable unions of \Borel partial transversals, it follows that there
  is a decreasing sequence $\sequence{E_n}[n \in \N]$ of \Borel equivalence relations such that every
  $E_n$ has only countably many classes, each of which is a partial transversal of $F_n$ of diameter
  at most $1 / n$. For each $n \in \N$, set $G_n = F_n \intersection G$, and note that if $\sequence{x_i}
  [i \le k]$ is a $G_n$-path for which $\sequence{\equivalenceclass{x_i}{E_n}}[i \le k]$ is a
  $(G_n)_{E_n}$-cycle, then $x_0 \mathrel{(E_n \intersection F_n)} x_k$, so $x_0 = x_k$, thus 
  $\sequence{x_i}[i \le k]$ is a $G_n$-cycle. It follows that each $G_n$ has faithful cycles on $E_n$, 
  so $\sequence{E_n, G_n}[n \in \N]$ is a partition stratification of $G$.
  
  Conversely, suppose that $\sequence{E_n, G_n}[n \in \N]$ is a partition stratification of $G$.
  By Theorem \ref{preliminaries:coloring}, there is a \Borel edge $\N$-coloring $c$ of $G$.
  Let $H_n$ denote the subgraph of $G_n$ consisting of all $\pair{x}{y} \in G_n \setminus E_n$
  for which $c(x, y)$ is minimal both among natural numbers of the form $c(x', y)$ where $x 
  \mathrel{E_n} x' \mathrel{G} y$, and those of the form $c(x, y')$ where $x \mathrel{G} 
  y' \mathrel {E_n} y$. Then $\sequence{H_n}[n \in \N]$ is an increasing sequence of \Borel 
  graphs whose union is $G$. By Theorem \ref{preliminaries:hypersmooth}, to see that 
  $\inducedequivalencerelation{G}$ is hyperfinite, we need only check that the relations
  $F_n = \inducedequivalencerelation{H_n}$ are smooth. By Proposition \ref
  {preliminaries:countablesmooth}, it is sufficient to show that for all $n \in \N$, every equivalence 
  class of $E_n$ is a partial transversal of $F_n$.
  
  Suppose, towards a contradiction, that $k \in \N$ is least for which there is an injective $H_n$-path
  $\sequence{x_i}[i \le k]$ beginning and ending at distinct $E_n$-related points. The definition of 
  partition stratification ensures that $\sequence{\equivalenceclass{x_i}{E_n}}[i \le k]$ is not a 
  $(G_n)_{E_n}$-cycle, so there exists $0 < i < k$ for which $x_{i-1} \mathrel{E_n} x_{i+1}$. Set 
  $x = x_{i - 1}$, $y = x_i$, and $x' = x_{i + 1}$. Then $c(x, y) \neq c(x', y)$, so the definition of $H_n$ 
  ensures that $\neg x \mathrel{H_n} y$ or $\neg x' \mathrel{H_n} y$, the desired contradiction.
\end{propositionproof}

We say that properties $P$ and $Q$ of \Borel equivalence relations \definedterm{coincide} 
below a given \Borel equivalence relation $F$ if the family of \Borel equivalence relations on
\Polish spaces which are \Borel reducible to $F$ and satisfy $P$ is the same as the family of
\Borel equivalence relations on \Polish spaces which are \Borel reducible to $F$ and satisfy $Q$.

\begin{proposition} \label{hyperfiniteoncountable:stratification}
  Suppose that $X$ is a \Polish space, $E$ is a \Borel equivalence relation on $X$, and 
  $G$ is a \Borel graphing of $E$ which admits a partition stratification. Then countability and
  hyperfiniteness coincide below $E$.
\end{proposition}

\begin{propositionproof}
  Fix a partition stratification $\sequence{E_n, G_n}[n \in \N]$ of $G$.
  
  \begin{lemma} \label{hyperfiniteoncountable:convexclosure}
    There are only countably many injective $G$-paths between any two points.
  \end{lemma}
  
  \begin{lemmaproof}  
    Suppose, towards a contradiction, that $k \in \N$ is the least natural number for which there exist
    $x, y \in X$ between which there are uncountably many injective $G$-paths $\sequence{z_i}[i \le k]$
    from $x$ to $y$. Then for $n \in \N$ sufficiently large, there are uncountably many injective 
    $G_n$-paths $\sequence{z_i}[i \le k]$ from $x$ to $y$ with the further property that 
    $\sequence{\equivalenceclass{z_i}{E_n}}[i \le k]$ is an injective $\inducedgraph{(G_n)}{E_n}$-path.
    Fix such an injective $G_n$-path $\sequence{z_i}[i \le k]$ from $x$ to $y$ for which there are 
    uncountably many injective $G_n$-paths $\sequence{z_i'}[i \le k]$ from $x$ to $y$ inducing the 
    same injective $\inducedgraph{(G_n)}{E_n}$-path as $\sequence{z_i}[i \le k]$. The minimality of 
    $k$ then ensures that there are uncountably many injective $G_n$-paths $\sequence{z_i'}[i \le k]$ 
    from $x$ to $y$ inducing the same injective $\inducedgraph{(G_n)}{E_n}$-path as $\sequence{z_i}
    [i \le k]$ but avoiding $\set{z_i}[0 < i < k]$. Then for $n' \in \N$ sufficiently large, there are 
    uncountably many injective $G_n$-paths $\sequence{z_i'}[i \le k]$ from $x$ to $y$ inducing the 
    same injective $\inducedgraph{(G_n)}{E_n}$-path as $\sequence{z_i}[i \le k]$ but for which the 
    corresponding injective $\inducedgraph{(G_n)}{E_{n'}}$-path avoids $\set{\equivalenceclass{z_i}
    {E_{n'}}}[0 < i < k]$. It follows that there is such an injective $G_n$-path $\sequence{z_i'}[i \le k]$
    from $x$ to $y$ for which there are uncountably many such injective $G_n$-paths $\sequence{z_i''}
    [i \le k]$ inducing the same injective $\inducedgraph{(G_n)}{E_{n'}}$-path as $\sequence{z_i'}
    [i \le k]$. By one more appeal to the minimality of $k$, there are uncountably many such injective 
    $G_n$-paths $\sequence{z_i''}[i \le k]$ inducing the same injective $\inducedgraph{(G_n)}
    {E_{n'}}$-path as $\sequence{z_i'}[i \le k]$ but avoiding $\set{z_i'}[0 < i < k]$. Fix any such injective 
    $G_n$-path $\sequence{z_i''}[i \le k]$, and observe that $\sequence{z_1', \ldots, z_k' = z_k, \ldots, 
    z_0 = z_0'', z_1''}$ is an injective $G_n$-path inducing a $\inducedgraph{(G_n)}{E_{n'}}$-cycle,
    contradicting the fact that $G_{n'}$ has faithful cycles on $E_{n'}$.
  \end{lemmaproof}

  By Theorem \ref{preliminaries:LusinNovikov} and Proposition \ref{preliminaries:hyperfinite}, it is 
  sufficient to show that $E$ is hyperfinite on every \Borel set $B \subseteq X$ on which it is countable. 
  Towards this end, let $C$ denote the \definedterm{convex closure} of $B$ with respect to $G$, that 
  is, the set of points lying along an injective $G$-path between two points of $B$. As Theorem \ref
  {preliminaries:LusinNovikov} and Lemma \ref{hyperfiniteoncountable:convexclosure} ensure that 
  $C$ is also a \Borel set on which $E$ is countable, Proposition \ref
  {hyperfiniteoncountable:countablestratification} implies that $E$ is hyperfinite on $C$, thus
  hyperfinite on $B$.
\end{propositionproof}

In the treeable case, we can say even more.

\begin{proposition} \label{hyperfiniteoncountable:essential}
  Suppose that $X$ is a \Polish space, $E$ is a treeable \Borel equivalence relation on $X$, and
  countability and hyperfiniteness coincide below $E$. Then essential countability and essential
  hyperfiniteness also coincide below $E$.
\end{proposition}

\begin{propositionproof}
  Suppose that $Y$ is a \Polish space and $F$ is an essentially countable \Borel equivalence relation
  on $Y$ which admits a \Borel reduction $\phi \from Y \to X$ to $E$. Fix a \Polish space
  $Y'$ and a countable \Borel equivalence relation $F'$ on $Y'$ for which there is a \Borel 
  reduction $\psi \from Y \to Y'$ of $F$ to $F'$. Then the set $R_0 = \set{\pair{x}{\psi(y)}}[\phi(y) 
  = x]$ \definedterm{induces a partial injection} of $X / E$ into $Y' / F'$, in the sense that 
  $x_1 \mathrel{E} x_2 \iff y_1' \mathrel{F'} y_2'$, for all $\pair{x_1}{y_1'}, \pair{x_2}{y_2'} \in R_0$.
  
  The \definedterm{product} of the equivalence relations $E$ and $F'$ is the relation on $X \times Y'$ 
  given by $\pair{x_1}{y_1'} \mathrel{(E \times F')} \pair{x_2}{y_2'} \iff (x_1 \mathrel{E} x_2 \mathand 
  y_1' \mathrel{F} y_2')$.
  
  \begin{lemma}
    There is an $(E \times F')$-invariant \Borel set $R \supseteq R_0$ inducing a partial injection of 
    $X / E$ into $Y' / F'$.
  \end{lemma}
  
  \begin{lemmaproof}
    As the property of inducing a partial injection of $X / E$ into $Y' / F'$ is $\Piclass[1]
    [1]$-on-$\Sigmaclass[1][1]$ and closed under $(E \times F')$-saturation, by repeatedly applying 
    Theorem \ref{preliminaries:reflection}, we obtain \Borel sets $R_{n+1} \supseteq \saturation{R_n}
    {E \times F'}$ inducing \Borel partial injections of $X / E$ into $Y' / F'$. Define $R = \union[n \in \N]
    [R_n]$.
  \end{lemmaproof}

  As $F'$ is countable, Theorem \ref{preliminaries:LusinNovikov} ensures that the set $C = \image
  {\projection{X}}{R}$ is \Borel, and that there is a \Borel uniformization $\pi \from C \to Y'$ of 
  $R$. As any such function is necessarily a reduction of $E$ to $F'$ on $C$, it follows that $E$ is 
  essentially countable on $C$. An application of Theorem \ref{preliminaries:treeable} then yields a 
  \Borel set $D \subseteq C$, whose $E$-saturation is $C$, on which $E$ is countable. As 
  countability and hyperfiniteness coincide below $E$, it follows that $E$ is hyperfinite on $D$, and 
  one more application of Theorem \ref{preliminaries:LusinNovikov} yields a \Borel reduction of
  $\restriction{E}{C}$ to $\restriction{E}{D}$, so $E$ is essentially hyperfinite on $C$, thus $F$ is 
  essentially hyperfinite.
\end{propositionproof}

\section{Ideals} \label{ideals}

We say that a family $\calK$ of subsets of $\N \times \N$ is \definedterm{determined by
cardinalities on vertical sections} if $A \in \calK \iff B \in \calK$, whenever $\forall n \in \N
\ \cardinality{\verticalsection{A}{n}} = \cardinality{\verticalsection{B}{n}}$. 

For each family $\scrN \subseteq \powerset{\N}$ of subsets of the natural numbers, we use
$\closure{\scrN}$ to denote the closure of $\scrN$ under finite unions, and we define
$\cardinalityideal{\scrN} = \union[N \in \closure{\scrN}][\cardinalityideal{N}]$, where
\begin{equation*}
  \cardinalityideal{N} = \set{A \subseteq \N \times \N}[\forall n \in \N \ (\cardinality{\verticalsection{A}
    {n}} = \aleph_0 \implies n \in N)].
\end{equation*}
Note that every such family is both determined by cardinalities on vertical sections and an
\definedterm{ideal}, in the sense that it is closed under containment and finite unions.

\begin{proposition} \label{ideals:continuous}
  Suppose that $\scrN \subseteq \powerset{\N}$ is a family of pairwise almost disjoint infinite subsets
  of $\N$. Then there is a continuous function $\pi \from \scrN \to \powerset{\N \times \N}$ \Wadge 
  reducing $\scrM$ to $\cardinalityideal{\scrM}$, for all $\scrM \subseteq \scrN$.
\end{proposition}

\begin{propositionproof}
  Define $\pi(N) = N \times \N$. Given $\scrM \subseteq \scrN$, observe first that if $M \in \scrM$, then
  $\pi(M) \in \cardinalityideal{M} \subseteq \cardinalityideal{\scrM}$. Conversely, if $N \in
  \scrN$ and $\pi(N) \in \cardinalityideal{\scrM}$, then there exist $n \in \N$ and $M_1, \ldots,
  M_n \in \scrM$ such that $\pi(N) \in \cardinalityideal{M_1 \union \cdots \union M_n}$, in which
  case $N \subseteq M_1 \union \cdots \union M_n$. As $\scrN$ consists of pairwise almost disjoint 
  infinite sets, it follows that $N = M_i$, for some $1 \le i \le n$, thus $N \in \scrM$.
\end{propositionproof}

A weak converse to this result is provided by the following.

\begin{proposition} \label{ideals:Borel}
  Suppose that $\scrN \subseteq \powerset{\N}$ is a \Borel family of pairwise almost disjoint infinite
  subsets of $\N$. Then $\cardinalityideal{\scrN}$ is also \Borel.
\end{proposition}

\begin{propositionproof}
  Clearly $\cardinalityideal{\scrN}$ is analytic, so by Theorem \ref{preliminaries:bianalytic}, it is 
  sufficient to show that it is co-analytic. But this follows from Theorem \ref{preliminaries:unicity} and
  the fact that a set $A$ is in $\cardinalityideal{\scrN}$ if and only if there exist $k \in \N$
  and a unique subfamily $\calF$ of $\scrN$ of cardinality $k$ for which there is a finite subset $F$ 
  of a set in $\closure{\scrN}$ such that $\forall n \in \N \ (\cardinality{\verticalsection{A}{n}} = \aleph_0 
  \implies n \in F \union \union[\calF])$.
\end{propositionproof}

\section{Trees} \label{trees}

Suppose that $t_{i, n} \in \Cantorspace[n]$, for all $i < 2$ and $n \in \N$. Associated with 
$\sequence{t_{i, n}}[i < 2, n \in \N]$ are the graphs $T_n$ on $\Cantorspace[n]$ obtained 
recursively by letting $T_{n+1}$ be the union of the graph $\set{\pair{s \concatenation 
\singletonsequence{i}}{t \concatenation \singletonsequence{i}}}[i < 2 \mathand \pair{s}{t} 
\in T_n]$ with the singleton edge $\set{\pair{t_{i, n} \concatenation \singletonsequence{i}}
{t_{1 - i, n} \concatenation \singletonsequence{1 - i}}}[i < 2]$, as well as the set $T = 
\set{\pair{\emptyset}{\emptyset}} \union \union[n \in \N][T_n]$.

A straightforward induction shows that each $T_n$ is a tree on $\Cantorspace[n]$. It follows that if
the set $T$ is closed under initial segments, then the set $\branches{T} = \set{\pair{x}{y} \in 
\Cantorspace \times \Cantorspace}[\forall n \in \N \ \pair{\restriction{x}{n}}{\restriction{y}{n}} \in T_n]$ of 
\definedterm{branches} through $T$ is an acyclic compact graph admitting a partition stratification, 
thus all of its \Borel subgraphs admit partition stratifications as well. As equivalence relations induced 
by acyclic \Borel graphs are themselves \Borel (by Theorems \ref{preliminaries:bianalytic} and 
\ref{preliminaries:unicity}), Propositions \ref{hyperfiniteoncountable:stratification} and 
\ref{hyperfiniteoncountable:essential} therefore imply that essential countability and essential 
hyperfiniteness coincide below equivalence relations generated by such subgraphs.

The \definedterm{support} of a sequence $c \in \Cantorspace$ is given by $\support{c} =
\preimage{c}{1}$.

\begin{proposition} \label{ideals:symmetricdifference}
  Suppose that $t_{i, n} \in \Cantorspace[n]$ for all $i < 2$ and $n \in \N$, the corresponding set
  $T$ is closed under initial segments, $\pair{a}{b}, \pair{c}{d} \in \branches{T}$, and $\set{a, b} \neq 
  \set{c, d}$. Then $(\support{a} \symmetricdifference \support{b}) \intersection (\support{c} 
  \symmetricdifference \support{d})$ is finite.
\end{proposition}

\begin{propositionproof}
  For each $n \in \N$, let $T_n$ denote the tree on $\Cantorspace[n]$ associated with $\sequence
  {t_{i, n}}[i < 2, n \in \N]$, and note that if $n \in \support{a} \symmetricdifference \support{b}$, then
  the pair of restrictions $\pair{\restriction{a}{(n+1)}}{\restriction{b}{(n+1)}}$ is in $T_{n+1}$, from
  which it follows that $\set{\restriction{a}{(n+1)}, \restriction{b}{(n+1)}} = \set{t_{i,n} \concatenation
  \singletonsequence{i}}[i < 2]$, and therefore that $\set{\restriction{a}{n}, \restriction{b}{n}} =
  \set{t_{i,n}}[i < 2]$. In particular, if $n \in \N$ is sufficiently large that $\set{\restriction{a}{n},
  \restriction{b}{n}} \neq \set{\restriction{c}{n}, \restriction{d}{n}}$, then $n$ is not in $(\support{a} 
  \symmetricdifference \support{b}) \intersection (\support{c} \symmetricdifference \support{d})$.
\end{propositionproof}

Given a family $\calI$ of subsets of $\N$, let $\inducedequivalencerelation{\calI}$ denote the binary
relation on $\Cantorspace$ given by $c \mathrel{\inducedequivalencerelation{\calI}} d \iff \support
{c} \symmetricdifference \support{d} \in \calI$.

\begin{proposition} \label{ideals:graphing}
  Suppose that $t_{i, n} \in \Cantorspace[n]$ for all $i < 2$ and $n \in \N$, the corresponding set
  $T$ is closed under initial segments, and $\calI$ is an ideal on $\N$ containing all finite subsets
  of $\N$. Then $\inducedequivalencerelation{\calI} \intersection \branches{T}$ is a graphing of 
  $\inducedequivalencerelation{\calI} \intersection \inducedequivalencerelation{\branches{T}}$.
\end{proposition}

\begin{propositionproof}
  It is sufficient to show that if $n \in \N$, $\sequence{c_i}[i \le n]$ is an injective $\branches
  {T}$-path, and $c_0 \mathrel{\inducedequivalencerelation{\calI}} c_n$, then $\sequence{c_i}[i \le 
  n]$ is an $(\inducedequivalencerelation{\calI} \intersection \branches{T})$-path. Towards this
  end, appeal to Proposition \ref{ideals:symmetricdifference} repeatedly to obtain a finite set $F
  \subseteq \N$ containing $(\support{c_i} \symmetricdifference \support{c_{i+1}}) \intersection 
  (\support{c_j} \symmetricdifference \support{c_{j+1}})$, for all $i < j < n$. Then $\support{c_0} 
  \symmetricdifference \support{c_n}$ and $\union[i < n][\support{c_i} \symmetricdifference 
  \support{c_{i+1}}]$ agree off of $F$, so $\support{c_i} \symmetricdifference \support{c_{i+1}}
  \subseteq F \union (\support{c_0} \symmetricdifference \support{c_n})$, and is therefore
  in $\calI$, for all $i < n$.
\end{propositionproof}

\section{Density} \label{density}

Suppose that $t_{i, n} \in \Cantorspace[n]$, for all $i < 2$ and $n \in \N$, and the corresponding set
$T$ is closed under initial segments. Then for each $n \in \N$, there exist $k \in \set{0, \ldots, n}$
and $s \in \Cantorspace[n - k]$ with $t_{i, n + 1} = t_{i, k} \concatenation \singletonsequence{i} 
\concatenation s$. Conversely, if $k_n \in \set{0, \ldots, n}$ and $s_n \in \Cantorspace[n - k_n]$ for
all $n \in \N$, then the set $T$ associated with $\sequence{t_{i,n}}[i < 2, n \in \N]$, where 
$t_{i, 0} = \emptystring$ and $t_{i, n + 1} = t_{i, k_n} \concatenation \singletonsequence{i} 
\concatenation s_n$ for $i < 2$ and $n \in \N$, is closed under initial segments. We say that
$\sequence{k_n, s_n}[n \in \N]$ is \definedterm{suitable} if $k_n \in \set{0, \ldots, n}$ and $s_n \in 
\Cantorspace[n - k_n]$, for all $n \in \N$.

Fix an injective enumeration $\sequence{i_n, j_n}[n \in \N]$ of $\N \times \N$ with $i_0 = 0$, and 
let $e$ denote its inverse. We say that $\sequence{k_n, s_n}[n \in \N]$ is \definedterm{dense} (with
respect to our fixed enumeration) if for all $i, k \in \N$ and all $s \in \Cantortree$, there exists $n \in
\N$ such that $i = i_{n+1}$, $k = k_n$, and $s \extendedby s_n$.

The \definedterm{push-forward} of a family $\calK$ of subsets of $\N \times \N$ through $e$ is given
by $\pushforward{e}{\calK} = \set{\image{e}{A}}[A \in \calK]$.

\begin{proposition} \label{density:complexity}
  Suppose that $\sequence{k_n, s_n}[n \in \N]$ is dense and suitable, $T$ is the corresponding
  set, $C \subseteq \Cantorspace$ is comeager, $\calK$ is a family of subsets of $\N \times
  \N$ which is determined by cardinalities on vertical sections and invariant under finite alterations
  of the leftmost column, and $\calI = \pushforward{e}{\calK}$. Then there is a \Wadge reduction
  $\pi \from \powerset{\N} \to \restriction{\branches{T}}{C}$ of $\calI$ to 
  $\inducedequivalencerelation{\calI}$.
\end{proposition}

\begin{propositionproof}
  Fix dense open sets $U_n \subseteq \Cantorspace$ whose intersection is contained in $C$, and let
  $t_{i, n}$ denote the sequences associated with $\sequence{k_n, s_n}[n \in \N]$. We will 
  recursively construct natural numbers $\ell_n > 0$, in addition to natural numbers $n_u < \ell_n$ 
  and sequences $t_u \in \Cantorspace[\ell_n - n_u - 1]$ for $u \in \Cantorspace[n]$, from which we 
  define $\pi_{i,n} \from \Cantorspace[n] \to \Cantorspace[\ell_n]$ by $\pi_{i, n}(u) = t_{i, n_u} 
  \concatenation \singletonsequence{i} \concatenation t_u$, for $i < 2$ and $u \in \Cantorspace[n]$. 
  We will ensure that the following conditions hold:
  \begin{enumerate}
    \item $\forall i, j < 2 \forall u \in \Cantorspace[n] \ \pi_{i, n}(u) \strictlyextendedby 
      \pi_{i, n + 1}(u \concatenation \singletonsequence{j})$.
    \item $\forall i < 2 \forall u \in \Cantorspace[n+1] \ \extensions{\pi_{i, n+1}(u)} 
      \subseteq U_n$.
    \item $\forall u \in \Cantorspace[n] \ \pair{\pi_{0,n}(u)}{\pi_{1,n}(u)} \in T$.
    \item $\forall u \in \Cantorspace[n] \ i_n = i_{n_{u \concatenation \singletonsequence{1}}}$.
    \item $\forall u \in \Cantorspace[n] \ \support{\pi_{0, n + 1}(u \concatenation \singletonsequence{0})} 
      \symmetricdifference \support{\pi_{1, n + 1}(u \concatenation \singletonsequence{0})}$ \\ 
        \hspace*{45pt} $= \support{\pi_{0, n}(u)} \symmetricdifference \support{\pi_{1, n}(u)}$.
    \item $\forall u \in \Cantorspace[n] \ \support{\pi_{0, n + 1}(u \concatenation \singletonsequence{1})} 
        \symmetricdifference \support{\pi_{1, n + 1}(u \concatenation \singletonsequence{1})}$ \\ 
          \hspace*{45pt} $= (\support{\pi_{0, n}(u)} \symmetricdifference \support{\pi_{1, n}(u)}) \union 
            \set{n_{u \concatenation \singletonsequence{1}}}$.
  \end{enumerate}
  We begin by setting $\ell_0 = 1$, $n_\emptyset = 0$, and $t_\emptystring = \emptystring$.
  
  Suppose now that $n \in \N$ and we have already found $\ell_n$, $n_u$, and $t_u$, for $u \in 
  \Cantorspace[n]$. For each $u \in \Cantorspace[n]$, set $n_{u \concatenation \singletonsequence
  {0}} = n_u$, and fix $t_u' \in \Cantortree$ with the property that $\extensions{\pi_{i, n}(u) 
  \concatenation t_u'} \subseteq U_n$ for all $i < 2$. By density, there exists $n_{u \concatenation 
  \singletonsequence{1}} > 0$ with $i_n = i_{n_{u \concatenation \singletonsequence{1}}}$, 
  $n_u = k_{n_{u \concatenation \singletonsequence{1}} - 1}$, and $t_u \concatenation t_u'
  \extendedby s_{n_{u \concatenation \singletonsequence{1}} - 1}$. Define $\ell_{n+1} = \max_{u \in 
  \Cantorspace[n]} n_{u \concatenation \singletonsequence{1}} + 1$, and for each $u \in \Cantorspace
  [n]$, fix an extension $t_{u \concatenation \singletonsequence{0}} \in \Cantorspace[\ell_{n+1} -
  n_u - 1]$ of $t_u \concatenation t_u'$, as well as $t_{u \concatenation \singletonsequence{1}} \in 
  \Cantorspace[\ell_{n+1} - n_{u \concatenation \singletonsequence{1}} - 1]$.
  
  By condition (1), the function $\pi \from \powerset{\N} \to \Cantorspace \times \Cantorspace$ which 
  is given by $\pi(A) = \pair{\pi_0(\characteristicfunction{A})}{\pi_1(\characteristicfunction{A})}$, where 
  $\pi_i(c) = \union[n \in \N][\pi_{i, n}(\restriction{c}{n})]$ for all $i < 2$, is well-defined and continuous. 
  Condition (2) then ensures that $\image{\pi_i}{\Cantorspace} \subseteq C$ for all $i < 2$, so 
  $\image{\pi}{\powerset{\N}} \subseteq C \times C$, thus condition (3) implies that $\image{\pi}
  {\powerset{\N}} \subseteq \restriction{\branches{T}}{C}$. And conditions (4) - (6) ensure that for all 
  $N \subseteq \N$, the cardinalities of all but the leftmost vertical sections of $\preimage{e}
  {\support{\pi_0(\characteristicfunction{N})} \symmetricdifference \support{\pi_1(\characteristicfunction
  {N})}}$ and $\preimage{e}{N}$ agree, whereas the cardinalities of their leftmost vertical sections differ 
  by at most one, thus $N \in \calI$ if and only if $\pi(N) \in \inducedequivalencerelation{\calI}$.
\end{propositionproof}

As a consequence, we obtain the following.

\begin{theorem} \label{density:main}
  The family of treeable \Borel equivalence relations below which essential hyperfiniteness and
  the inexistence of a continuous embedding of $\Eone$ coincide has cofinal essential complexity.
\end{theorem}

\begin{theoremproof}
  By appealing to Proposition \ref{preliminaries:almostdisjoint}, we obtain a continuous injection $\pi 
  \from \Cantorspace \to \powerset{\N}$ into a family of pairwise almost disjoint infinite sets. Given a 
  \Borel set $B \subseteq \Cantorspace$, set $\scrN = \image{\pi}{B}$. Proposition \ref
  {ideals:continuous} then ensures that $B$ is \Wadge reducible to $\cardinalityideal{\scrN}$, and 
  Proposition \ref{ideals:Borel} implies that the latter is \Borel, thus the same holds of the ideal
  $\calI = \pushforward{e}{\cardinalityideal{\scrN}}$.
  
  Fix a dense suitable sequence $\sequence{k_n, s_n}[n \in \N]$, and let $T$ be the
  associated set. Proposition \ref{ideals:graphing} then ensures that the equivalence 
  relation $E$ on $\Cantorspace$ given by $E = \inducedequivalencerelation{\calI} \intersection 
  \inducedequivalencerelation{\branches{T}}$ is generated by a \Borel subgraph of an
  acyclic compact graph, and since it is clearly \Borel, by Theorem \ref{preliminaries:dichotomy} 
  we need only check that its essential complexity is at least the complexity of $B$.
  
  Towards this end, suppose that $Y$ is a \Polish space, and $F$ is a \Borel equivalence relation on
  $Y$ for which there is a \Borel reduction $\phi \from X \to Y$ of $E$ to $F$. Proposition \ref
  {preliminaries:continuous} then yields a comeager \Borel set $C \subseteq X$ on which $\phi$ is 
  continuous. By Proposition \ref{density:complexity}, there is a \Wadge reduction $\psi \from
  \powerset{\N} \to \restriction{\branches{T}}{C}$ of $\calI$ to $E$, and it follows that $(\phi \times
  \phi) \composition \psi$ is a \Wadge reduction of $\calI$ to $F$, so $B$ is \Wadge reducible to
  $F$, thus the essential complexity of $E$ is at least the complexity of $B$.
\end{theoremproof}

\section{Anti-basis results} \label{antibasis}

We say that a family $\calF$ of \Borel equivalence relations on \Polish spaces is \definedterm
{unbounded} if for every \Borel equivalence relation $E$ on a \Polish space, there is a \Borel equivalence relation $F \in \calF$ which is not \Borel reducible to $E$.

We say that the non-linearity of \Borel reducibility is \definedterm{captured} off of a family $\calF$ of
\Borel equivalence relations if every non-essentially-hyperfinite \Borel equivalence relation is
incompatible with a \Borel equivalence relation outside of $\calF$.

\begin{theorem} \label{antibasis:nonlinearity}
  Suppose that $\calF$ is a class of \Borel equivalence relations whose complement contains
  unboundedly many \Borel equivalence relations below which essential hyperfiniteness and
  the inexistence of a \Borel reduction of $\Eone$ coincide, as well as unboundedly many 
  \Borel equivalence relations to which $\Eone$ does not admit a \Borel reduction. Then the 
  non-linearity of \Borel reducibility is captured off of $\calF$.
\end{theorem}

\begin{theoremproof}
  Suppose that $E$ is a \Borel equivalence relation compatible with every \Borel equivalence relation
  outside of $\calF$. Fix a \Borel equivalence relation $E'$, outside of $\calF$ and not
  \Borel reducible to $E$, below which essential hyperfiniteness and the inexistence of a
  \Borel reduction of $\Eone$ coincide. In addition, fix a \Borel equivalence relation $E''$, outside 
  of $\calF$ and not \Borel reducible to $E$, to which $\Eone$ does not admit a \Borel reduction. As 
  $E$ is \Borel reducible to $E''$, it follows that $\Eone$ is not \Borel reducible to $E$. As $E$ is 
  \Borel reducible to $E'$, it therefore follows that $E$ is essentially hyperfinite.
\end{theoremproof}

The following corollary strengthens \cite[Theorem 2]{KechrisLouveau}.

\begin{theorem} \label{antibasis:incompatible}
  Suppose that $X$ is a \Polish space and $E$ is a \Borel equivalence relation on $X$.
  Then the following are equivalent:
  \begin{enumerate}
    \item The relation $E$ is not essentially hyperfinite.
    \item There is a \Borel equivalence relation on a \Polish space which is incompatible with $E$.
    \item The family of \Borel equivalence relations on \Polish spaces which are incompatible with $E$ 
      under \Borel reducibility has cofinal essential complexity.
  \end{enumerate}
\end{theorem}

\begin{theoremproof}
  As $(3) \implies (2)$ is trivial and $(2) \implies (1)$ is a consequence of Theorem \ref
  {preliminaries:Ezero}, it is sufficient to show $(1) \implies (3)$. Towards this end, note that
  Theorem \ref{density:main} ensures that the family of \Borel equivalence relations 
  on \Polish spaces below which essential hyperfiniteness and the inexistence of a \Borel 
  reduction of $\Eone$ coincide has cofinal essential complexity, and Theorems \ref
  {preliminaries:groupaction} and \ref{preliminaries:nonreducible} imply that the family of \Borel 
  equivalence relations on \Polish spaces to which $\Eone$ does not \Borel reduce has cofinal 
  essential complexity. So, by Theorem \ref{antibasis:nonlinearity}, the family of \Borel equivalence 
  relations on \Polish spaces which are incompatible with $E$ under \Borel reducibility has cofinal 
  essential complexity.
\end{theoremproof}

We can now establish our primary result.

\begin{theorem}
  Suppose that $\calF$ is a dichotomical class of \Borel equivalence relations on \Polish spaces of 
  bounded essential complexity. Then every equivalence relation in $\calF$ is smooth.
\end{theorem}

\begin{theoremproof}
  Fix a \Borel equivalence relation $F$ witnessing that $\calF$ is dichotomical. Then $F$ is
  necessarily essentially hyperfinite, since otherwise Theorem \ref{antibasis:incompatible} would
  yield a \Borel equivalence relation outside of $\calF$ and incompatible with $F$. But Theorem
  \ref{preliminaries:Ezero} then implies that every relation in $\calF$ is smooth.
\end{theoremproof}

\begin{acknowledgments}
  We would like to express our gratitude to Renaud Bonnafoux and Jean-Marc Tomczyk, 
  who took over Dominique Lecomte's teaching while he spent a semester in M\"{u}nster, during
  which the bulk of this work was completed. We would also like to thank Alain Louveau and
  Gabriel Debs for their remarks on our work.
\end{acknowledgments}

\vfill\eject

\bibliographystyle{amsalpha}
\bibliography{bibliography}

\end{document}